\documentclass[12pt]{amsart}
\usepackage{amscd}
\usepackage{amssymb}
\usepackage{amstext}
\usepackage{amsthm}
\usepackage{mathrsfs}
\usepackage{latexsym}
\usepackage{xcolor}
\numberwithin{equation}{section}
\usepackage{hyperref}

\newcommand{\cc}{\mathbb{C}}
\newcommand{\N}{\mathbb{N}}

\newcommand{\D}{\mathbb{D}}

\newcommand{\Bn}{\mathbb B_n}
\newcommand{\Sn}{\mathbb S_n}
\newcommand\set[1]{\left\{#1\right\}}
\providecommand{\abs}[1]{\lvert#1\rvert}
\providecommand{\Abs}[1]{\Bigl\lvert#1\Bigr\rvert}

\def\A2w{A^2_\rho}
\def\F2w{\mathcal F^2_\psi}

\providecommand{\normw}[1]{\lVert#1\rVert_\rho}

\newcommand{\less}{\lesssim}
\newcommand{\gess}{\gtrsim}
\newcommand{\asym}{\asymp}
\newcommand{\vh}[1]{\langle #1\rangle}
\providecommand{\norm}[1]{\lVert#1\rVert}

\def\B{\mathcal B}
\def\d{\partial}

\def\Dhat{\widehat{\mathcal{D}}}

\newtheorem{Thm}{Theorem}[section]
\newtheorem{theorem}[Thm]{Theorem}
\newtheorem{lemma}[Thm]{Lemma}
\newtheorem{proposition}[Thm]{Proposition}

\newtheorem*{thmA}{Theorem A}
\newtheorem*{thmB}{Theorem B}

\theoremstyle{definition}
\newtheorem{remark}[Thm]{Remark}

\begin{document}
\title[On the Bergman projections acting on $L^\infty$]{On the Bergman projections acting on $L^\infty$ in the unit ball $\Bn$}
\author{Van An Le}


\address{Aix--Marseille University, CNRS, Centrale Marseille, I2M, Marseille, France}
\address{University of Quynhon, Department of Mathematics, 170 An Duong Vuong, Quy Nhon, Vietnam}
\email{vanandkkh@gmail.com}


\begin{abstract} Given a weight function, we define the Bergman type projection with values in the corresponding weighted Bergman space on the unit ball $\Bn$ of $\cc^n, n>1$. We characterize the radial weights such that this projection is bounded from $L^\infty$ to the Bloch space $\B$.
\end{abstract}

\maketitle

\section{Introduction and main result}
Let $\cc^n$ denote the $n$-dimensional complex Euclidean space. For any two points
$z=(z_1, \ldots, z_2),\, w=(w_1, \ldots, w_n)$ in $\cc^n$, we use the well-known notation 
$$\vh{z,w}=z_1\overline{w_1}+\cdots+z_n\overline{w_n}\quad \text{ and }\quad
\abs z=\sqrt{\vh{z,z}}.$$ Let $\Bn=\{z\in \cc^n:|z|<1\}$ be the unit ball, and let 
$\Sn=\{z\in \cc^n: \abs z=1\}$ be the unit sphere in $\cc^n$.
Denote by $H(\Bn)$ the space of all holomorphic functions on the
unit ball $\Bn$.
Let $dv$ be the normalized volume measure on $\Bn$. The normalized surface measure on $\Sn$ will be denoted by $d\sigma$.

Let $\rho$ be a positive and integrable function on $[0,1)$. We  extend it to $\Bn$ by $\rho(z)=\rho(\abs z)$, and call such $\rho$ a radial weight function. The weighted Bergman space $\A2w$ is the space of functions $f \text{ in } H(\Bn)$ such that 
 $$ \normw{f}^2 =\int_{\Bn}\abs{f(z)}^2\rho( z)dv(z)<\infty. $$

Let $\rho$ be a radial weight and $X$ be a space of measurable functions on $\Bn$. The Bergman type projection $P_\rho$ acting on $X$ is given by
$$P_\rho f(z) = \int_{\Bn} K_\rho(z,w)f(w)\rho(w)dv(w), \qquad z\in \Bn, f\in X,$$
where $K_\rho(z,w)$ is the reproducing kernel of the weighted Bergman space $\A2w$.

  When $\rho$ is the standard radial weight $\rho(z)=(1-\abs z^2)^\alpha, \alpha>-1$, the corresponding projection is denoted by $P_\alpha$.
  
  A radial weight $\rho$ belongs to the class $\Dhat$ if $\widehat{\rho}(r)\lesssim \widehat\rho(\frac{1+r}{2})$ for all $r\in [0,1)$, where $\widehat \rho(r)=\int_r^1 \rho(s)ds$.
 
 The study of small Bergman spaces in higher dimensions began in 2018 in our work \cite{An}. Projections play a crucial role in studying operator theory on spaces of analytic functions. Bounded
analytic projections can also be used to establish duality relations and to obtain useful equivalent norms in spaces of analytic functions. Hence the boundedness of projections is an
interesting topic which has been studied by many authors in recent years \cite{C-P,Dostanic1,Dostanic,pelaez20161,pelaez2019}. In \cite{pelaez20161}, Pel\'aez and R\"atty\"a  considered the projection $P_{\rho_1}$ acting on $L^p_{\rho_2}(\D), 1\le p<\infty$, when  two weights $\rho_1, \rho_2 $ are in the class $\mathcal R$ of so called regular weights. A radial weight $\rho$ is regular if $\widehat{\rho}(r)\asymp (1-r)\rho(r), r\in (0,1)$. Recently, in 2019, they extended these results to the case where $\rho_1\in \Dhat$, $\rho_2$ is radial \cite{pelaez2019}.

In this text, we are going to study the projections acting on the space $L^\infty$.
Let us recall that the Bloch space of $\Bn$, denoted by $\B(\Bn)$, or simply by $\B$, is the space of holomorphic functions $f$ in $\Bn$ such that
$$\sup _{z\in \Bn}(1-\abs z^2)\abs{Rf(z)}<\infty,$$
where $$Rf(z)=\sum_{j=1}^n z_j\frac{\d f}{\d z_j}(z)$$ is the radial derivative of $f$ at $z\in \Bn$. 
 In the one dimensional case, the Bloch space consists of analytic functions $f$ on $\D$ such that
 $$\sup _{z\in \D}(1-\abs z^2)\abs{f'(z)}<\infty,$$
 and is denoted by $\B(\D)$.

 In the case of standard radial weight, we have the following result.
\begin{thmA}
For any $\alpha>-1$, the Bergman type projection $P_\alpha$ is a
bounded linear operator from $L^\infty$ onto the Bloch space $\B$.
\end{thmA}
 See \cite[Theorem 5.2]{zhudisk} for the proof in the case of one variable and \cite[Theorem 3.4]{zhuball} for the proof in the case of several variables.
 
In \cite{pelaez2019}, Pel\'aez and R\"atty\"a obtained an interesting result in the one dimensional case.
\begin{thmB}
Let $\rho$ be a radial weight. Then the projection $P_\rho: L^\infty(\D)\to \B(\D)$ is bounded if and only if $\rho\in \Dhat$.
\end{thmB}
We extend this theorem to the case of several variables and obtain the following result.
\begin{theorem}\label{thm6.1.3}
Let $\rho$ be a radial weight. Then the projection $P_\rho: L^\infty\to \mathcal B$ is bounded if and only if $\rho\in \Dhat$.
\end{theorem}


Throughout this text, the notation $U(z)\lesssim V(z) $ (or equivalently $V(z)\gess U(z)$) means that there is a positive constant $C$ such that $U(z)\le CV(z)$ holds for all $z$ in the set in question, which may be a space of functions or a set of numbers. If both $U(z)\less V(z)$ and $V(z)\less U(z)$, then we write $U(z)\asym V(z)$.

\section{Some auxiliary lemmas}
To prove Theorem \ref{thm6.1.3} we need several auxiliary lemmas.
\begin{lemma}\label{lem5}
Let $\rho$ be a radial weight. Then the following conditions are equivalent:
\begin{itemize}
\item[(i)]$\rho\in \Dhat$;
\item[(ii)] There exist $C=C(\rho)>0$ and $\beta_0=\beta_0(\rho)>0$ such that
$$\widehat{\rho}(r)\le C\left (\dfrac{1-r}{1-t}\right )^\beta \widehat{\rho}(t),\qquad 0\le r\le t<1,$$
for all $\beta\ge \beta_0$;
\item[(iii)] The asymptotic equality 
$$\int_0^1s^x \rho(s)ds\asymp \widehat\rho\left (1-\frac1x\right ), \qquad x\in [1,\infty),$$
is valid;
\item[(iv)] There exist  $C_0=C_0(\rho)>0$ and $C=C(\rho)>0$ such that
$$\widehat\rho(0)\le C_0\widehat\rho(\frac12)$$
and $\rho_n\le C\rho_{2n}$ for all $n\in \N$.
\end{itemize}
\end{lemma}
This lemma can be found in \cite{pelaez2016lem}.

\begin{lemma}\label{lem3}
If $$f(z)=\sum_{n=0}^\infty a_jz^j \in H^p, \qquad 0<p\le 2,$$
then 
$$\sum_{j=0}^{\infty}(j+1)^{p-2}\abs{a_j}^p\less\norm f_p^p.$$
\end{lemma}
\begin{lemma}
Let $\set{a_j}$ be a sequence of complex numbers such that $\sum j^{q-2}\abs{a_j}^q<\infty$ for some $q, 2\le q <\infty$. Then the function $f(z)=\sum_{n=0}^\infty a_jz^j$ is in $H^q$, and
$$\norm f_q^q\less\sum_{j=0}^{\infty}(j+1)^{q-2}\abs{a_j}^q.$$
\end{lemma}

Two above lemmas are the classical Hardy-Littlewood inequalities, which can be found, for example, in Duren's book \cite[Theorem 6.2 and 6.3]{Duren}.

\begin{lemma}
Let $\rho$ be a radial weight. Then the reproducing kernel $K_\rho(z,w)$ is given by
$$K_\rho(z,w)=\frac12\sum_{d=0}^\infty\dfrac{(d+n-1)!}{d!n!\rho_{2n-1+2d}}\vh{z,w}^d, \qquad z, w\in \Bn,$$
where 
$$\rho_{x}=\int_0^1 t^{x}\rho(t)dt,\qquad x\ge 1.$$
\end{lemma}
\begin{proof}
By the multinomial formula (see \cite[(1.1)]{zhuball}), we have that
$$\vh{z,w}^d=\sum_{\beta\in \N^n, \abs \beta =d}\frac{d!}{\beta!}z^\beta\bar w^{\beta}, \qquad z,w \in \cc^n.$$
Hence, for $\alpha\in \N^n, \abs \alpha=d$,
$$\int_{\Sn}\xi^\alpha\vh{z,\xi}^dd\sigma(\xi)= \sum_{\beta\in \N^n, \abs \beta =d}\frac{d!z^\beta}{\beta!}\int_{\Sn}\xi^\alpha\bar \xi^\beta d\sigma(\xi), \quad z\in \Bn.$$
By Lemma 1.11 in \cite{zhuball},
$$\int_{\Sn}\xi^\alpha\bar \xi^\beta d\sigma(\xi)=
\begin{cases}
0	& \text{ if } \alpha\ne \beta,\\
\dfrac{\alpha!(n-1)!}{(d+n-1)!}& \text{ if } \alpha= \beta,
\end{cases}$$
and we obtain
\begin{align*}
\int_{\Sn}\xi^\alpha\vh{z,\xi}^dd\sigma(\xi)
	&=\dfrac{d!}{\alpha!} z^\alpha\int_{\Sn}\xi^\alpha\bar \xi^\alpha d\sigma(\xi)\\
	&=\dfrac{d!}{\alpha!}\dfrac{\alpha!(n-1)!}{(d+n-1)!}z^\alpha\\
	&=\dfrac{d!(n-1)!}{(d+n-1)!}z^\alpha, \quad z\in \Bn.
\end{align*}
Therefore, for $\alpha\in \N^n, \abs \alpha =d$ we have
\begin{align*}
\int_{\Bn}w^\alpha\vh{z,w}^d\rho(w)dv(w)
	&=2n\int_0^1 t^{2n-1+2d}\rho(t)dt\int_{\Sn}\xi^\alpha\vh{z,\xi}^dd\sigma(\xi)\\
	&=\dfrac{2d!n!\rho_{2n-1+2d}}{(d+n-1)!}z^\alpha,\qquad z \in \Bn,
\end{align*} 
It follows that
\begin{equation}\label{star}
z^{\alpha}=\dfrac{(d+n-1)!}{2d!n!\rho_{2n-1+2d}}\int_{\Bn}w^\alpha\vh{z,w}^d\rho(w)dv(w),\qquad z\in \Bn.
\end{equation}

Since $\rho(t)>0, 0<t<1$, we have $\rho_s\ge C_\varepsilon (1-\varepsilon)^s$ for every $\varepsilon>0$. Given $z\in \Bn$, we have
\begin{align*}
&\int_{\Bn}\Abs{\frac12\sum_{d=0}^\infty\dfrac{(d+n-1)!}{d!n!\rho_{2n-1+2d}}\vh{z,w}^d}^2\rho(w)dv(w)\\
	&=\frac14\sum_{d_1,d_2\ge 0}\dfrac{(d_1+n-1)!(d_2+n-1)!}{d_1!d_2!(n!)^2\rho_{2n-1+2d_1}\rho_{2n-1+2d_2}}\int_{\Bn}\vh{z,w}^{d_1} \vh{w,z}^{d_2}\rho(w)dv(w)\\
	&=\frac14\sum_{d_1,d_2\ge 0}\dfrac{(d_1+n-1)!(d_2+n-1)!}{d_1!d_2!(n!)^2\rho_{2n-1+2d_1}\rho_{2n-1+2d_2}} \times \\
	& \hspace*{5cm}\times\int_{\Bn}\sum_{\abs \beta=d_2}w^\beta\bar z^\beta\frac{d_2!}{\beta!}\vh{z,w}^{d_1}\rho(w) dv(w)\\
	&=\frac12\sum_{d\ge 0}\left (\dfrac{(d+n-1)!}{d!n!}\right )\dfrac{1}{\rho_{2n-1+2d}^2}\sum_{\abs \beta=d}\dfrac{(d!)^2}{\beta!}\dfrac{n!\rho_{2n-1+2d}}{(d+n-1)!}z^\beta\bar z^\beta\\
	&=\frac12\sum_{d\ge 0}\dfrac{(d+n-1)!}{n!\rho_{2n-1+2d}}\sum_{\abs \beta=d}\dfrac{z^\beta\bar z^\beta}{\beta!}=\frac12\sum_{d\ge 0}\dfrac{(d+n-1)!}{d!n!\rho_{2n-1+2d}}\abs z^{2d}<\infty.
\end{align*}
Thus, the function $w\mapsto \frac12\sum_{d= 0}^\infty\dfrac{(d+n-1)!}{d!n!\rho_{2n-1+2d}}\vh{w,z}^d$ belongs to $\A2w$.

By \eqref{star} and by continuity, for every $f\in A^2_\rho(\Bn)$, 
$$
f(z)=\int_{\Bn}f(w)\left (\frac12\sum_{d=0}^\infty\dfrac{(d+n-1)!}{d!n!\rho_{2n-1+2d}}\vh{z,w}^d\right )\rho(w)dv(w),\qquad z\in \Bn,$$
which implies our conclusion. 
\end{proof}

\section{Proof of main result}
It suffices to consider only the case $n>1$.

\begin{proposition}\label{pro1}
If $\rho\in \Dhat$, then the projection $P_\rho: L^\infty\to \mathcal B$ is bounded, where $P_\rho$ is defined by
$$P_\rho\varphi(z)=\int_{\Bn}K_\rho(z,w)\varphi(w)\rho(w)dv(w), \qquad \varphi\in L^\infty, z\in \Bn.$$
\end{proposition}

\begin{proof}
We have 
$$K_\rho(z,w)=\frac12\sum_{d=0}^\infty\dfrac{(d+n-1)!}{d!n!\rho_{2n-1+2d}}\vh{z,w}^d.$$
Hence, for a fixed $w\in \Bn$,
\begin{align*}
RK_\rho(z,w)
	&=\sum_{j=1}^{n}z_j\dfrac{\d K_\rho(z,w)}{\d z_j}\\
	&=\sum_{j=1}^{n}z_j\frac{\d}{\d z_j}\left (\frac12\sum_{d=0}^\infty\dfrac{(d+n-1)!}{d!n!\rho_{2n-1+2d}}\vh{z,w}^d\right )\\
	&=\frac12\sum_{j=1}^{n}z_j\sum_{d=0}^\infty\dfrac{(d+n-1)!}{d!n!\rho_{2n-1+2d}}d \bar w_j\vh{z,w}^{d-1}\\
	&=\frac12\sum_{d=1}^\infty\dfrac{(d+n-1)!}{(d-1)!n!\rho_{2n-1+2d}}\vh{z,w}^d\\
	&=\frac1{2}\sum_{d=1}^\infty \dfrac{\Gamma(d+n)}{\Gamma(d)\Gamma(n+1)\rho_{2n-1+2d}}\vh{z,w}^d.
\end{align*}
Now, given $\varphi \in L^\infty$, let
$$f(z):=P_\rho\varphi(z)=\int_{\Bn}K_\rho(z,w)\varphi(w)\rho(w)dv(w), \qquad z\in \Bn.$$
For all $z\in \Bn$ we have
\begin{align}
\abs{Rf(z)}
	&=\Abs{\int_{\Bn}RK_\rho(z,w)\varphi(w)\rho(w)dv(w)}\notag\\
	&\le \int_{\Bn}\abs{RK_\rho(z,w)}\abs{\varphi(w)}\rho(w)dv(w)\notag\\
	&\le \norm\varphi_\infty \int_{\Bn} \abs{RK_\rho(z,w)}\rho(w)dv(w).\label{5.2}
\end{align}

Set $$g(\lambda)=\sum_{d=1}^\infty\dfrac{\Gamma(d+n)}{\Gamma(d)} \dfrac{\lambda^{d-1}}{\rho_{2n-1+2d}}, \qquad \lambda \in \D.$$
Since $\rho(t)>0, 0<t<1,$ $g$ is analytic in the unit disc.
Then 
\begin{equation}\label{5.3}
RK_\rho(z,w)= \dfrac{\vh{z,w}}{2\Gamma(n+1)}g(\vh{z,w}).
\end{equation}
Next we consider the reproducing kernel $K^1_{\rho}(z,w)$ of the Bergman space in the unit disc with the weight $\rho$. We have
$$K^1_\rho(z,w)=\dfrac12\sum_{d=0}^\infty \dfrac{(z\overline{w})^{d}}{\rho_{2d+1}}.$$ 
Furthermore,
\begin{align*}
\dfrac{\d^n}{\d z^n}K^1_\rho(z,w)
	&=\dfrac12\sum_{d=n}^\infty \dfrac{\Gamma(d+1)(z\overline w)^{d-n}\overline w^n}{\Gamma(d-n+1)\rho_{2d+1}}\\
	&=\dfrac12\sum_{s=1}^\infty\dfrac{\Gamma(s+1)}{\Gamma(s)}\dfrac{(z\overline w)^{s-1}\overline w^n}{\rho_{2s+2n-1}}\\
	&=\dfrac12 g(z\overline w) \overline w^n.
\end{align*}
By a result of Pel\'aez and R\"atty\"a (\cite[Theorem 1 (ii)]{pelaez20161}), we have
$$\int_{\D}\Abs{\dfrac{\d^n}{\d z^n}K^1_\rho(z,w)}(1-\abs z^2)^{n-2} dA(z)\asymp \int_0^{\abs w}\dfrac{dt}{\widehat\rho(t)(1-t)^2}, \qquad \dfrac{1}{2}\le \abs w<1,$$
where $\widehat\rho(t)=\int_t^1\rho(s)ds$.

Thus, 
$$\int_{\D}\abs{g(z\overline{w})}(1-\abs{z}^2)^{n-2}dA(z)\asymp \int_0^{\abs w}\dfrac{dt}{\widehat\rho(t)(1-t)^2}, \qquad \dfrac{1}{2}\le \abs w<1.$$
Since $g$ is analytic in the unit disc, we have
\begin{equation}\label{5.4}
\int_{\D}\abs{g(z\overline{w})}(1-\abs{z}^2)^{n-2}dA(z)\lesssim 1+ \int_0^{\abs w}\dfrac{dt}{\widehat\rho(t)(1-t)^2}, \qquad w\in \D.
\end{equation}
Now, by \eqref{5.3}, we have
\begin{align*}
\int_{\Bn}\abs{RK_\rho(z,w)}\rho(w)dv(w)
	&\lesssim \int_{\Bn}\abs{g(\vh{z,w})}\rho(w)dv(w)\\
	&\asymp\int_0^1 r^{2n-1}\rho(r)\left (\int_{\Sn}\abs{g(\vh{rz,\xi})}d\sigma(\xi)\right )dr.	
\end{align*}
By \cite[Lemma 1.9]{zhuball} and the unitary invariance of $d\sigma$, we have
$$\int_{\Sn}\abs{g(\vh{rz,\xi})}d\sigma(\xi)\asymp \int_{\D}\abs{g(r\abs z \lambda)}(1-\abs \lambda^2)^{n-2}dA(\lambda).$$
Thus, by \eqref{5.4} we obtain
\begin{align*}
\int_{\Bn}\abs{RK_\rho(z,w)}&\rho(w)dv(w)\\
	&\lesssim \int_0^1 r^{2n-1}\rho(r)\left (1+ \int_0^{r\abs z}\dfrac{dt}{\widehat\rho(t)(1-t)^2}\right )dr\\
	&\lesssim 1+\int_0^{\abs z}\dfrac{1}{\widehat{\rho}(t)(1-t)^2}\left (\int_{t/\abs z}^1 r^{2n-1}\rho(r)dr\right )dt\\
	&\lesssim 1+\int_0^{\abs z}\dfrac{\widehat{\rho}(t/\abs z)}{\widehat{\rho}(t)}\dfrac{dt}{(1-t)^2}\lesssim \dfrac{1}{1-\abs z}, \quad z\in \Bn.
\end{align*}

By \eqref{5.2} we obtain now that
$$
\abs{Rf(z)}	\lesssim \norm{\varphi}_\infty \dfrac{1}{1-\abs z^2},\quad z\in \Bn,
$$
and, hence, 
$$\sup _{z\in \Bn}(1-\abs z^2)\abs{Rf(z)}\lesssim \norm\varphi_\infty. $$

It is easy to see that $$\abs{f(0)}\lesssim \norm\varphi_\infty.$$
Therefore, $P_\rho$ is bounded. The Proposition~\ref{pro1} is proved.
\end{proof}



\begin{proposition}\label{pro2}
Suppose that the projection $P_\rho: L^\infty\to \mathcal B$ is bounded. Then $\rho\in \Dhat$.
\end{proposition}

\begin{proof}
Given $\xi\in \Sn$ and $w\in \Bn$, let us consider a function $g$ given by
$$g(\lambda)=RK_\rho(\lambda \xi,w), \qquad \lambda\in \mathbb D.$$
Then 
$$g(\lambda)=\sum_{d=1}^{\infty}c_d\vh{\xi,w}^d\lambda^d,$$
where $c_d=\dfrac1{2n}\dfrac{\Gamma(d+n)}{\Gamma(d)\Gamma(n)\rho_{2n-1+2d}}$.
By the Hardy--Littlewood inequality (see Lemma~\ref{lem3}) we have 
$$\sum_{d=1}^{\infty}\dfrac{c_d\abs{\vh{\xi,w}}^d}{d+1}\lesssim \int_0^{2\pi}\abs{g(e^{i\theta})}\dfrac{d\theta}{2\pi}=\int_0^{2\pi}\abs{RK_\rho(e^{i\theta}\xi,w)}\dfrac{d\theta}{2\pi}.$$
Integrating both sides of the above inequality over $\xi\in \Sn$ we obtain
\begin{align*}
\sum_{d=1}^{\infty}\dfrac{c_d}{d+1}\int_{\Sn}\abs{\vh{\xi,w}}^d\,d\sigma(\xi)
	&\lesssim \int_{\Sn}\int_0^{2\pi}\abs{RK_\rho(e^{i\theta}\xi,w)}\dfrac{d\theta}{2\pi}d\sigma(\xi)\\
	&=\int_{\Sn}\abs{RK_\rho(\xi,w)}\,d\sigma(\xi).
\end{align*}
By the unitary invariance of $d\sigma$ and \cite[Lemma 1.9]{zhuball}, we have
\begin{align*}
\int_{\Sn}\abs{\vh{\xi,w}}^d\,d\sigma(\xi)
	&=\abs w^d\int_{\Sn}\abs{\xi_1}^d\,d\sigma(\xi)\\
	&=(n-1)\abs w^d\int_{\mathbb D}(1-\abs z^2)^{n-2}\abs z^d\,dA(z)\\
	&=(n-1)\pi\abs w^d \int_0^1(1-t)^{n-2}t^{d/2}dt\\
	&\asymp\dfrac{\Gamma(\frac d2+1)\Gamma(n)}{\Gamma(\frac d2+n)}\abs w^d.
\end{align*}
Hence,
\begin{align*}
\int_{\Sn}\abs{RK_\rho(\xi,w)}\,d\sigma(\xi)
	&\gtrsim \sum_{d=1}^{\infty}\dfrac{c_d}{d+1}\dfrac{\Gamma(\frac d2+1)\Gamma(n)}{\Gamma(\frac d2+n)}\abs w^d\\
	&=\frac1{2n}\sum_{d=1}^{\infty}\dfrac{\Gamma(d+n)\Gamma(\frac d2+1)}{(d+1)\Gamma(d)\Gamma(\frac d2+n)\rho_{2n-1+2d}}\abs w^d.
\end{align*}
Since 
$$\dfrac{\Gamma(d+n)\Gamma(\frac d2+1)}{(d+1)\Gamma(d)\Gamma(\frac d2+n)}\asymp 1,$$
we get
$$\int_{\Sn}\abs{RK_\rho(\xi,w)}\,d\sigma(\xi)\gtrsim \frac1{2n}\sum_{d=1}^{\infty} \dfrac{\abs w^d}{\rho_{2n-1+2d}}, \qquad w\in \Bn.$$
Therefore, for $z\in \Bn$, we have
\begin{align*}
\int_{\Bn} \abs{RK_\rho(z,w)}\rho(w)dv(w)
	&= 2n\int_0^1 r^{2n-1}\rho(r)\int_{\Sn}\abs{RK_\rho(z,r\xi)}\,d\sigma(\xi) \,dr\\
	&=2n\int_0^1 r^{2n-1}\rho(r)\int_{\Sn}\abs{RK_\rho(\xi,rz)}\,d\sigma(\xi) \,dr\\
	&\gtrsim \sum_{d=1}^{\infty} \dfrac{\abs z^d}{\rho_{2n-1+2d}}\int_0^1 r^{2n-1+d}\rho(r) dr\\
	&=\sum_{d=1}^{\infty} \dfrac{\rho_{2n-1+d}}{\rho_{2n-1+2d}}\abs z^d.
\end{align*}
Thus, 
\begin{align*}
\sup_{z\in Bn}(1-\abs z^2)
	&\int_{\Bn} \abs{RK_\rho(z,w)}\rho(w)dv(w)\\
	&\gtrsim \sup_{z\in \Bn}(1-\abs z)\sum_{d=1}^{\infty} \dfrac{\rho_{d+2n-1}}{\rho_{2d+2n-1}}\abs z^d\\
	&\ge \sup _{N\in \N} \dfrac 1N \sum_{d=1}^{N}\dfrac{\rho_{d+2n-1}}{\rho_{2d+2n-1}}\left (1-\dfrac 1N\right )^d\\
	&\gtrsim \sup _{N\in \N} \dfrac 1N \sum_{d=1}^{N}\dfrac{\rho_{d+2n-1}}{\rho_{2d+2n-1}}.
\end{align*}
Since $P_\rho$ is bounded, 
$$\sup_{z\in Bn}(1-\abs z^2)\int_{\Bn} \abs{RK_\rho(z,w)}\rho(w)dv(w)<\infty.$$
Given $N\ge 2n$, we obtain that 
$$1\gtrsim \dfrac 1{4N-2n} \sum_{d=3N-n+1}^{4N-2n}\dfrac{\rho_{d+2n-1}}{\rho_{2d+2n-1}}\ge \dfrac{1}{4N}(N-n) \dfrac{\rho_{4N}}{\rho_{6N}},$$
and, hence,
$$\rho_{6N}\gtrsim \rho_{4N}.$$
If $8N\le k< 8N+8, N\ge 2n+8$, then
$$\rho_k\le\rho_{8N}\lesssim \rho_{12N}\lesssim \rho_{18N}\le \rho_{2k},$$
and by Lemma \autoref{lem5} we conclude that $\rho\in \Dhat$.
\end{proof}

From Propositions \ref{pro1} and \ref{pro2}, we obtain the conclusion of  Theorem \ref{thm6.1.3}.

\begin{remark}
The method given herein combined with our results in \cite{An} can be used to generalize to the unit ball case the $L^p$ estimates proved in \cite{pelaez20161} in the unit disk case. This will be the object of a forthcoming paper.
\end{remark}

\vspace{0.2cm}
\noindent \textbf{Acknowledgments.} I am deeply grateful to Alexander Borichev and El Hassan Youssfi for their help and many suggestions during the preparation of this paper.


\end{document}